\documentclass[1pt,notitlepage,twoside,a4paper]{amsart}

\usepackage{amsmath,amssymb,enumerate}

\usepackage{epsfig,fancyhdr,color}

\usepackage{amssymb}
\usepackage{amsmath,amsthm}    
\usepackage{latexsym}
\usepackage{amscd} 
\usepackage{psfrag}
\usepackage{graphicx} 
\usepackage[latin1]{inputenc}  
\usepackage[all]{xy} 
\usepackage[mathcal]{eucal}

\theoremstyle{definition}

\theoremstyle{remark}

\def\interieur#1{\mathord{\mathop{\kern 0pt #1}\limits^\circ}}

\definecolor{NoteColor}{rgb}{1,0,0}

\title[A Commentary on Teichm\"uller's paper]{A Commentary on Teichm\"uller's paper \emph{\"Uber Extremalprobleme der konformen Geometrie (On extremal problems in conformal geometry)}}\author[N. A'Campo and A. Papadopoulos]{Norbert A'Campo
and Athanase Papadopoulos}

\address{N. A'Campo: Universit\"at Basel,  Mathematisches Institut, 
\\
Spiegelgasse 1, 4051 Basel, Switzerland
\\
 email:\,\tt{norbert.acampo@gmail.com}}
\address{A. Papadopoulos: Institut de Recherche Math\'ematique Avanc\'ee, UMR 7501
\\
Universit{\'e} de Strasbourg and CNRS,\\
7 rue Ren\'e Descartes, 67084 Strasbourg Cedex, France,
\\
And: Max Planck Institute for Mathematics,
\\
Vivatsgasse 7, 53111 Bonn
\\
email:\,\tt{papadop@math.unistra.fr}}

 \date{\today}

\begin{document}

 \begin{abstract}

We comment on the paper \emph{\"Uber Extremalprobleme der konformen Geometrie (On extremal problems in conformal geometry)} \cite{T23} by Teichm\"uller, published in 1941. This paper contains ideas on a wide generalization of his previous work on the solution of extremal problems in conformal geometry. The generalization concerns at the same time the fields of function theory,  topology and algebra.

 \medskip


 \medskip

 \medskip

The final version of this paper will appear as a chapter in Volume VI of {\it the Handbook of Teichm\"uller theory}. This volume is dedicated to the memory of Alexander Grothendieck.

\end{abstract}

\bigskip

\maketitle

\tableofcontents

We comment on the paper \emph{\"Uber Extremalprobleme der konformen Geometrie (On extremal problems in conformal geometry)} \cite{T23} by Teichm\"uller, published in 1941. This paper may be considered as a collection of thoughts that generalize some of the ideas expressed in Teichm\"uller's paper \cite{T20} (see also the commentary in \cite{T20C}) where he uses quadratic differentials to solve an extremal problem about quasiconformal mappings, which is a substantial extension of the problem on extremal quasiconformal mappings between rectangles solved by Gr\"otzsch. In the present paper, \cite{T23},  Teichm\"uller expresses the idea that quadratic differentials may be used to provide solutions to other extremal problems in conformal geometry. This idea involves the introduction of a new structure at distinguished points of a surface.  An ``order'' for series expansions of the local coordinate charts is fixed at a distinguished point. In this setting, if a point is distinguished without further requirements, then the associated quadratic differential has a simple pole at that point. But sometimes,  at the distinguished point, the extremal problem requires a function which has a fixed value for its $n$ first derivatives. Then, the quadratic differential has a pole of order $n+1$ at that point. In the geometric language introduced later on by Ehresmann, the jet structure of a certain order at a distinguished point is fixed.  At the end of this commentary, we shall mention some of Strebel's works that further develop this point of view.

The paper \cite{T23}  is difficult to read, because the ideas are formulated vagely. This makes the reviewing of this article a laborious task. A note by the editors of the journal where the paper was published,  on the first page of the paper, informs the reader that this paper is ``obviously unfinished,'' and warns him that  ``unreasonably high demands are made on the reader's cooperation and imagination,'' that ``the assertions are not even stated precisely with rigour, neither proofs nor even any clues are given,'' that some things ``not of fundamental importance [...] occupy
a broad space for something almost unintelligible, while far too scarce hints are provided for
 fundamentally important individual examples.'' The editors declare that 
if they decided to publish the paper, ``despite all lack that distinguishes the work against the other papers in this journal [Deutsche Mathematik],
it is to bring up for discussion the thoughts contained therein relating to the theme of estimates for schlicht functions.''
In turn, the author declares at the beginning of the paper: ``Because I only have a limited vacation time at my disposal, I cannot give reasons for many things, but only assert them.''

Teichm\"uller's aim is to show how some of his ideas expressed in his previous papers -- some of them in the papers on moduli, and others in papers on algebra -- are related to each other, that they may lead to  general concepts, and that they are applicable to various situations, in particular to the coefficient problem for univalent functions.  
 
 Teichm\"uller starts by commenting on the fact that function theory is closely related to topology and algebra. For instance, one is led, in dealing with function-theoretic questions, to prove new generalizations of the Riemann-Roch theorem. He describes a situation where one needs for that Lie theory, rather than differential geometry. He then explains his choice of the notion of ``principal domain''\index{principal region} (Hauptbereich), used in a  previous paper \cite{T20}, to denote a Riemann surface with distinguished points. This notion is given a more general meaning here, and the stress now is on the marked points rather than on the supporting surface. The author makes an analogy with a situation in algebra, that has already considered in \cite{T4}. To say things more precisely, one is given three objects, $\mathfrak{A}, \mathfrak{A}_1,\mathfrak{A}_2$. In the geometric case,  $\mathfrak{A}$ is a principal domain,  $\mathfrak{A}_1$ is the support of the principal domain, that is, the underlying Riemann surface (with no distinguished points) and $\mathfrak{A}_2$ the set of distinguished points. In the algebraic context,  $\mathfrak{A}$ is a normal (Galois) extension of a field,  $\mathfrak{A}_1$ is a cyclic  field extension and $\mathfrak{A}_2$ is the set of generators of the Galois group. Teichm\"uller considers the following situation: At a distinguished point on a Riemann surface, one chooses a local coordinate $z$, and assumes that the other local coordinates $\tilde{z}$ are of the form
  \[\tilde{z}=z+a_{m+1}z^{m+1}+a_{m+2}z^{m+2}+\ldots
\]
  where $m$ is an integer. This implies (by using the chain rule) that, at such a point, the first $m$ derivatives of a map, computed in local coordinates, are well defined, that is, they do not depend on the choice of the admissible local coordinates. Such a distinguished point is said to be of \emph{order} $m$. 
 Using a modern language, the author is fixing, at a finite number of distinguished points, the order-$m$ jets of maps at  each such point. Teichm\"uller further states that it is practical to think of such a point of oder $m$ as $m+1$ points which are infinitely close together. (Recall that in order to compute the  $m$-th derivative, one may consider certain quotients of values of the function at $m+1$ points and then takes limits.)
 
 There is an analogous description of distinguished points  at  the boundary of the surface. The author mentions such a description without going into any specific details. He then explains how surfaces with distinguished points transform under appropriate maps. (Recall that holomorphic functions are differentiable, therefore they act on jets.) In fact, Teichm\"uller considers the case of maps of a collection of disjoint subdomains (with distinguished points) of a Riemann surface to a  collection of disjoint subdomains (with distinguished points) of another Riemann surface. The maps under consideration are subject to some topological requirements. The reader may think of these subdomains as obtained through a decomposition of the  surface defined by the trajectory structure of a quadratic differential. Teichm\"uller mentions applications  to the Bieberbach coefficient problem. He mentions some explicit ones, and he refers to his previous papers \cite{T200} (see also the commentary \cite{T200C}).

 Later on, Lie theory is used. At each point, there is a continuous infinite-dimensional group $\mathfrak{G}_0$ of local parameter transformations, and a filtration
  \[\mathfrak{G}_0\supset \mathfrak{G}_1\supset \mathfrak{G}_2\supset \ldots
  \]
  where $\mathfrak{G}_m$ is a normal subgroup of $\mathfrak{G}_0$ whose elements are transformations $A$ of the form
  
   \[\tilde{z}=Az= z+a_{m+1}z^{m+1}+a_{m+2}z^{m+2}+\ldots
\]

Then, Lie algebras are defined.

The notion of ``topological determination'' (an expression which the author assigns  in \cite{T20} to an object we call today a marking) is given a much more general meaning. A generalized version of the Riemann-Roch theorem, which  takes into account the new notion of distinguished points, including points on the boundary, is also given. The formula involves sums over these points of the relative dimensions of some quotients of Lie algebras.
A notion of ``extremal problem in conformal geometry'' is formulated in such a wider setting. 

Then there is a long discussion on conformal embeddings of annuli (which he calls ring domains) on general surfaces (``higher principal domains'').\footnote{The reader familiar with the theory of pseudo-holomorphic curves, and especially with Floer homology, will recognize one of the basic ideas of this theory (in dimension two), where the embeddings of annuli control all the situation. In the theory of pseudo-holomorphic curves the study of embeddings of cylinders, in particular those that interpolate the Lagrangian submanifolds, are of fundamental importance.}
Teichm\"uller believes that the idea of studying embeddings of cylinders is hidden in the papers of Gr\"otzsch.\index{Gr\"otzsch, H.}

Teichm\"uller notes that the results apply to abstract function fields instead of Riemann surfaces, and that estimates on the coefficients of a univalent function may be obtained through a method involving extremal mappings associated with quadratic differentials with some prescribed poles.  He states that the proof is more complicated than the usual proof of the area problem, but that it has the advantage of admitting generalizations to higher degree coefficient problems. Teichm\"uller is alluding here to the famous Bieberbach conjecture, which was one of his main objects of interest. We state the problem and some developments below. In the same paper, he makes relations with several classical problems, including the question of finding the Koebe domain of a family of holomorphic functions defined on the disk, that is, the largest domain contained in the image of every function in the family, and the so-called ``Faber trick.'' Faber, in a paper  published in 1922 \cite{Faber}, found an early application of the length-area method\index{Gr\"otzsch-Ahlfors method}\index{length-area method} to the question of boundary correspondence of conformal mappings.

A translation of the paper \cite{T200} in which the length-area\index{Gr\"otzsch-Ahlfors method}\index{length-area method} method and geometric methods, such as properties of conformal moduli, are applied to the study of quasiconformal mappings and Riemann surfaces, and, in which, the length area method is used, will appear in Volume VII of this Handbook. (See also the commentary \cite{T200C}.) The length-area method is also used by Teichm\"uller in the papers \cite{T20}, \cite{T29} and \cite{T31}.

A few words are now in order for the Bieberbach conjecture, also known as the  Coefficient Problem for schlicht functions.\footnote{The German word ``schlicht'' is sometimes used in the English literature, and, of course, it is the word used by Teichm\"uller.}  These are univalent  (that is, holomorphic and injective)  functions defined on the unit disc $D=\{z\in \mathbb{C} \vert \ \vert z\vert <1\}$ by a Taylor series expansion: 
 \[f(z)=\sum_{n=0}^\infty a_nz^n
  \]
 normalized by $a_0=0$ and $a_1=1$. The Bieberbach conjecture, formulated by Bieberbach in 1916 and proved fully by Louis de Branges in 1984, says that the coefficients of such a series satisfy the inequalities
 \[\vert a_n\vert \leq n\]
 for any $n\geq 2$.
 In fact, Bieberbach proved in his paper \cite{BB} the case $n=2$, that is, he proved  $\vert a_2\vert \leq 2$. He also showed that equality is attained for the functions of the form $K_\theta (z)=z/(1-e^{i\theta}z)^2$ for $\theta\in\mathbb{R}$. This is the so-called ``Koebe function''  
 \[k(z)=z/(1-z)^2=z+2z^2+3z^3+\ldots
 \]
 composed with a rotation. In a footnote (\cite{BB} p. 946), Bieberbach suggested, without apparent motivation, that the value $n=\vert a_n(K_\theta)\vert$ might be an upper bound for all the functions satisfying the given assumptions. This is how the Bieberbach conjecture originated. Before its final proof by de Branges, a long list of good mathematicians worked on this problem, including L\"owner, Nevanlinna, Goluzin, Grunsky, Littlewood-Paley, Milin, Garabedian-Schiffer, Pederson-Ozawa,  Pederson-Schifferand and Jenkins.  Important  advances were made and results towards solving the conjecture were obtained. In particular, L\"owner proved in 1923 the second (and difficult) step that $\vert a_3\vert \leq 3$, Garabedian and Schiffer proved in 1955 that $\vert a_4\vert\leq 4$, Pederson and Ozawa proved in 1968 that $\vert a_6\vert\leq 6$, and Pederson and Schiffer proved in 1972 that $\vert a_5\vert\leq 5$. De Branges' proof, published in 1984, was a result of a new approach that uses operator theory. The proof was simplified, later on, and operator theory removed.
 
 Bieberbach's result is related to the so-called ``area theorem,'' which is referred to by Teichm\"uller in the present paper, a  theorem which gives the so-called ``Koebe quarter theorem,'' saying that the image of any univalent function $f$ from the unit disc of $\mathbb{C}$ onto a subset of $\mathbb{C}$ contains the disc of center $f\left( 0\right)$ and radius $\vert f^{\prime}\left( 0 \right)\vert /4$. In this geometric form of the Bieberbach conjecture, Teichm\"uller introduced the techniques of extremal quasiconformal mappings and quadratic differentials. This is the content of the last part of this paper.
 
 Now we come to some developments of Teichm\"uller's theory. In a series of works, Jenkins (cf. \cite{J1}, \cite{J2},\cite{J3}, \cite{J5} and others) developed an approach to the coefficient theorem that uses quadratic differentials and which is based on the works of Teichm\"uller and Gr\"otzsch.  In a 1962 ICM talk \cite{J3}, Jenkins writes: ``Teichm\"uller enunciated the intuitive principle that the solution of a certain type of extremal problem for univalent functions is determined by a quadratic differential for which the following prescriptions hold. If the competing mappings are to have a certain fixed point the quadratic differential will have a simple pole there. If in addition fixed values are required for the first $n$ derivatives of competing functions the quadratic differential will have a pole of order  $n+1$ at the point. He proved a coefficient result which represents a quite special case of the principle but did not obtain any general result of this type. The General Coefficient Theorem was presented originally as an explicit embodiment of Teichm\"uller's principle, that is, the competing functions were subjected to the normalizations implied by the above statement.'' 
In \cite{J1}, Jenkins has a similar quote concerning Teichm\"uller's result from \cite{T14}:  ``Teichm\"uller enunciated the principle that the solution of a certain type of extremal problem in geometric function theory is in general associated with a quadratic differential. [...] Teichm\"uller was led to this principle by abstraction from the numerous results of Gr\"otzsch [...] and by his considerations on quasiconformal mappings \cite{T20}. He applied this principle in certain concrete cases, the most important of which was his coefficient theorem \cite{T14} which is the most penetrating explicit result known in the general coefficient problem for univalent functions.''

In \cite{J3}, Jenkins states a theorem, which is in the spirit of  the results stated by Teichm\"uller in the paper which we review here, but in a precise form. The result concerns a Riemann surface of finite type equipped with a quadratic differential with a  decomposition of the surface into subdomains defined by the trajectory structure of this differential. There is a mapping from each of these subdomains onto non-overlapping subdomains of the surface, and these mappings are subject to conditions on preservation of poles, certain coefficients, and also to some topological conditions. The conclusion of the theorem is then an inequality that involves the coefficients of the quadratic differential at poles of order greater than one and those of the mapping, with a condition for the inequality to be an equality. This condition states that the function must be an isometry for the metric induced on the surface by the quadratic differential. It is followed by a detailed analysis of the equality case.  Jenkins also refers to numerous specific applications of such a result and this is in line with what Teichm\"uller has envisioned. (See \cite{J1}, p. 278-279.)

In conclusion, let us note for the readers of this Handbook that Thurston also thought about the Coefficient Problem, with an approach involving  tools which are familiar to  Teichm\"uller theorists: M\"obius transformations, the Schwarzian derivative, quasiconformal mappings and the universal Teichm\"uller space, cf. \cite{Th}.
 
 \medskip
 
 \noindent {\bf Acknowledgements} The authors thank M. Brakalova-Trevithick for her suggestions on this commentary. The second author acknowledges the support of the Max-Plank-Institute for Mathematics (Bonn) where this work was done.

\end{document}